\input amstex 
\documentstyle{amsppt}
\magnification1100
\NoRunningHeads
\refstyle{A}

\let\u=\underline

\define\1{{\bold 1}}
\define\C{{\Bbb C}}
\define\Z{{\Bbb Z}}

\define\wt{\operatorname{wt}}
\define\sh{\operatorname{sh}}
\define\lt{\operatorname{\ell \!\text{{\it t\,}}}}
\define\lspan{\operatorname{\Bbb F\text{-span}}}

\def\sq{\lower.3ex\vbox{\hrule\hbox{\vrule height1.2ex depth1.2ex\kern2.4ex 
\vrule}\hrule}\,}


\topmatter
\title{A basis of the basic
$\frak{sl}(3,\C)\sptilde$-module}
\endtitle
\author Arne Meurman and Mirko Primc
\endauthor
\address
Univ\. of Lund,
Dept\. of Mathematics, 
S-22100 Lund,
Sweden
\endaddress
\email
arnem\@maths.lth.se
\endemail
\address
Univ\. of Zagreb,
Dept\. of Mathematics,
Bijeni\v{c}ka 30, Zagreb,
Croatia
\endaddress
\email
primc\@math.hr
\endemail
\thanks Partially supported by a grant {}from the G\"oran
Gustafsson
Foundation for Research in Natural Sciences and Medicine and
by the Ministry of Science 
and Technology of the Republic of Croatia, grant 037002.
  \endthanks
\thanks This paper was circulated in 1997 and posted on
QA/9812029.
  \endthanks
\abstract
J.~Lepowsky and R.~L.~Wilson initiated the approach to 
combinatorial Rogers-Ramanujan type identities via the
vertex operator constructions of representations of affine
Lie algebras. In this approach the first new combinatorial 
identities were discovered by
S. Capparelli through the construction of the level 3 standard
$A^{(2)}_2$-modules. We obtained several 
infinite series of new combinatorial identities through 
the construction of all standard $A^{(1)}_1$-modules; the
identities associated to the fundamental modules coincide
with the two Capparelli identities.
In this paper we extend our construction to  the basic
$A^{(1)}_2$-module and, by using the principal specialization
of the Weyl-Kac character formula, we obtain a Rogers-Ramanujan type 
combinatorial identity for colored partitions. 
The new combinatorial identity indicates the next level of complexity
which one should expect in Lepowsky-Wilson's approach
for affine Lie algebras of higher ranks, say for $A^{(1)}_n$, $n\geq 2$,
in a way parallel to the next level of complexity seen when passing
{}from the Rogers-Ramanujan identities (for modulus $5$) to
the Gordon identities for odd moduli $\geq 7$.

\endabstract

\endtopmatter


\document

\head{ Introduction}\endhead

J.~Lepowsky and R.~L.~Wilson gave in \cite{LW} a Lie-theoretic 
interpretation and proof of the classical Rogers-Rama\-nu\-jan 
identities in terms of representations of the affine Lie algebra 
$\tilde\goth g=\goth {sl}(2,\Bbb C)\sptilde$. The identities are 
obtained by expressing in two ways the principal characters of 
vacuum spaces  for the principal Heisenberg subalgebra of 
$\tilde\goth g$. The product sides follow {}from the principally 
specialized Weyl-Kac character formula; the sum sides follow {}from 
the vertex operator construction of bases parametrized by 
partitions satisfying difference 2 conditions. Very roughly 
speaking, for a level 3 standard $\tilde\goth g$-module $L$
with a highest weight vector $v_0$, Lepowsky and Wilson construct 
$\Cal Z$-operators $\{Z(j)\mid j\in\Bbb Z\}$ which commute with 
the action of the principal Heisenberg subalgebra, and show that
$$
Z(j_1)Z(j_2)\dots Z(j_s)v_0,\qquad j_1\leq j_2\leq\dots\leq j_s <0
$$ 
is a spanning set of the vacuum space $\Omega_L$ of the principal
Heisenberg subalgebra. This spanning set is reduced to a basis by 
using the vertex operator formula
$$
\align
&\sum_{i\geq 0}\, a_i (Z(m-i)Z(n+i)+Z(n-i)Z(m+i))\\
&=-\tfrac38 (-1)^m\delta_{m+n,0}\pm (-1)^{\tfrac12+m+n}Z(m+n)\\
\endalign
$$ 
(the coefficients $a_i$ are defined by the binomial expansion
$(1+z)\sp{1/3}(1-z)\sp{-1/3}=\sum_{i\geq 0}\, a_i z\sp{i}$),
allowing them to erase {}from the spanning set monomials containing
the ``leading terms'' $Z(r)Z(r)$ and $Z(r-1)Z(r)$ of these relations. 
What is left is a spanning set of vectors $Z(j_1)Z(j_2)\dots Z(j_s)v_0$
which contain no factor of the form
$$
Z(j)Z(j),\quad Z(j-1)Z(j), \qquad j<0,
$$
or, equivalently, which satisfy the difference $2$ conditions
$j_{p+1}-j_p\geq 2$ (cf. \cite{A}). Finally, Lepowsky and Wilson 
prove the linear independence of this set of vectors, and this gives a
Lie-theoretic proof of the Rogers-Ramanujan identities.

Lepowsky-Wilson's approach is also possible for other affine Lie algebras
and for other constructions of vertex operators, as in
 \cite{C}, \cite{LP}, \cite{Ma} and  \cite{Mi}, for example.
In this approach the first new combinatorial identities were discovered by
S. Capparelli  through the construction of level 3 standard
$A^{(2)}_2$-modules in the principal picture. We obtained in \cite{MP2}
several 
infinite series of new combinatorial identities through 
the homogeneous construction of all standard $A^{(1)}_1$-modules; for the
fundamental modules and the $(1,2)$-specialization the identities 
coincide with the two Capparelli identities.
In this paper we follow the ideas developed in \cite{MP1} and \cite{MP2}
and  construct a basis of the basic
$A^{(1)}_2$-module parametrized by colored
partitions.

In order to
describe our main result, let $\goth g =\goth{sl}(3,\Bbb C)$ 
and let $e_i$, $h_i$, $f_i$ be the Chevalley
generators of $\goth g$. Set $\goth h=\Bbb C\text{-span}\, \{h_1,h_2\}$ and
let 
$B$ be the ordered basis 
$$
X_1=[e_1,e_2], X_2=e_1, X_3=e_2, X_4=h_1, X_5=h_2, X_6=f_2, X_7=f_1,
X_8=[f_2,f_1],
$$
$$
X_1\succ X_2\succ X_3\succ X_4\succ X_5\succ X_6\succ X_7\succ X_8.
$$
Let $\tilde\goth g$ be the affine Lie 
algebra associated with $\goth g$ (cf. \cite{K}), spanned by
elements $x(n)$, $x\in \goth g$, $n\in \Bbb Z$, the canonical central element
$c$ and
a derivation $d$. Then we have a Poincar\'e-Birkhoff-Witt spanning set
$$
X_{i_1}(j_1)X_{i_2}(j_2)\dots X_{i_s}(j_s)v_0,\quad
 j_1\leq j_2\leq\dots\leq j_s <0, \ 
X_{i_r}\preccurlyeq X_{i_{r+1}} \text{ if } j_r=j_{r+1},
$$ 
of the basic $\tilde\goth g$-module $L(\Lambda_0)$ 
with a highest weight vector $v_0$. We may refer to monomials
$X_{i_1}(j_1)X_{i_2}(j_2)\dots X_{i_s}(j_s)$ of the above form as
ordered monomials in the universal enveloping algebra $U(\tilde\goth g)$.

The above Poincar\'e-Birkhoff-Witt spanning set can be reduced to a
basis in the following way:
Let $\Cal R$ be the set of ordered monomials of the form
$$\align
&X_{i_1}(j)X_{i_2}(j), \text{ with ``colors'' } i_1i_2:
 11,21,22,31,32,33,42,43,44,51,52,53,54,55,\\
&\qquad \qquad\qquad\qquad\qquad \qquad\qquad\qquad
62,64,65,66,73,74,75,76,77,85,86,87,88,\\
&X_{i_1}(j-1)X_{i_2}(j), \text{ with ``colors'' } i_1i_2:
  11,12,13,14,15,16,17,18,22,24,26,27,28,\\
&\qquad \qquad\qquad\qquad\qquad \qquad\qquad\quad
33,35,36,37,38,47,48,56,58,66,68,77,78,88,\\
&X_3(j-1)X_4(j)X_1(j), \quad X_8(j-1)X_4(j-1)X_6(j);
\text{ \ for all } j<0.
\endalign$$

We shall say that an ordered monomial $u=X_{i_1}(j_1)X_{i_2}(j_2)\dots
X_{i_s}(j_s)$
 satisfies the difference $\Cal R$ conditions if for any given $d\in \Cal R$ 
the monomial $u$ does not contain as factors all the simple
factors of $d$, or, by abuse of language, if $u$ does not contain $d$ 
as a factor.  For example, $u=X_1(-2)X_5(-1)X_3(-1)$ does not
satisfy the difference $\Cal R$ conditions since it contains both of the
simple
factors of $d=X_1(-2)X_3(-1)\in\Cal R$.
On the other hand, $X_1(-6)X_5(-3)X_3(-1)$ obviously satisfies the
difference $\Cal R$ conditions since the ``differences'' between the ``parts''
$X_1(-6)$, $X_5(-3)$ and $X_3(-1)$ are ``big enough'', and the monomial
$X_1(-6)X_5(-3)X_3(-1)$ cannot contain any $d$ {}from $\Cal R$.

Now we can state our main result:
\proclaim{Theorem A}
 The set of vectors
$$
X_{i_1}(j_1)X_{i_2}(j_2)\dots X_{i_s}(j_s)v_0\in L(\Lambda_0),
$$
where the monomials  $X_{i_1}(j_1)X_{i_2}(j_2)\dots X_{i_s}(j_s)$ are ordered 
and satisfy the difference $\Cal R$ conditions, is a basis of the basic
$\tilde\goth g$-module $L(\Lambda_0)$.
\endproclaim

By analogy with the Rogers-Ramanujan case, we start with 
the vertex operator formula $X_1(z)^2|L(\Lambda_0)=0$ (cf. \cite{LP}),
and extract the coefficients of the powers of $z$:
$$
\sum_{i+j=n}\, X_1(i)X_1(j)= 0\quad\text{ on }\quad L(\Lambda_0).
$$ 
Using the adjoint action of $\goth g$
on $X_1(z)^2$ we get a $27$-dimensional space of relations
which annihilate $L(\Lambda_0)$. The set 
of quadratic monomials 
listed above is the set of leading terms of the components
of these relations (Lemma 1). We also construct two 
relations with the leading terms $X_3(j-1)X_4(j)X_1(j)$ and
$X_8(j-1)X_4(j-1)X_6(j)$ (Lemma 8).
 By using  these relations we can reduce the
Poincar\'e-Birkhoff-Witt spanning set to a basis,
 the proof is given at the end of the next section.

The proof of linear independence is based on  the same ideas that we used in  
the $\goth g=\frak{sl}(2,\C)$ case, the main result being Theorem 11 which is
analogous to Theorem 9.1 in \cite{MP2}. The main ingredient of the
proof is a construction of $(64+35+35+27)$-dimensional space of 
relations among relations (Proposition 3), which is obtained {}from the vertex
operator algebra
structure on $L(\Lambda_0)$ (cf. \cite{FLM} and \cite{MP2, Section 8}). 
These relations among relations are first reformulated in terms of
embeddings (Proposition 10, in a way the analogue of Lemma 9.2 in \cite{MP2}),
the general statement of Theorem 11 is obtained by using  Lemma 9.4 in
\cite{MP2}.
In the next section we introduce all necessary notions and facts needed in
the proof, except that we for a few details refer to \cite{MP2}.

As a consequence of Theorem A, and the principally specialized 
Weyl-Kac character formula, we obtain a combinatorial identity
for colored partitions of Rogers-Ramanujan type (Theorem B).

The new combinatorial identity indicates the next level of complexity
which one should expect in Lepowsky-Wilson's approach
for affine Lie algebras of higher ranks, say for $A^{(1)}_n$, $n\geq 2$,
in a way parallel to the next level of complexity seen when passing
{}from the Rogers-Ramanujan identities (for modulus $5$) to
the Gordon identities for odd moduli $\geq 7$.
This is due to the appearance of two additional cubic terms in the set
$\Cal R$ of
``mainly difference 2 conditions'' --- a new combinatorial phenomenon not seen
in the $A^{(1)}_1$ case. These cubic terms make the proof
of linear independence ``as complicated as'' the proof in the case
of level $2$ standard $A^{(1)}_1$-modules, suggesting, at least as far as
the complexity of partitions is concerned, some kind of duality between
``level $1$ rank $2$ identities'' and ``level $2$ rank $1$ identities''.

It is clear that Lepowsky-Wilson's approach for higher rank 
affine Lie algebras should lead to combinatorial Rogers-Ramanujan 
type identities for colored partitions. It should be noted that similar
combinatorial identities also appear as a consequence of the classical
$q$-series approach, see, for example, \cite{AAG} and the references
therein. 


\head{A basis of the basic module}\endhead

As was already stated,
let $\goth g =\goth{sl}(3,\Bbb F)$ and let $e_i$, $h_i$, $f_i$ be the Chevalley
generators of $\goth g$ and $\goth h=\lspan \{h_1,h_2\}$,
where $\Bbb F=\Bbb C$ or any other field of characteristic 0. Let $B$ be the
ordered basis 
$$\gather
X_1=[e_1,e_2], X_2=e_1, X_3=e_2, X_4=h_1, X_5=h_2, X_6=f_2, X_7=f_1,
X_8=[f_2,f_1],\\
X_1\succ X_2\succ X_3\succ X_4\succ X_5\succ X_6\succ X_7\succ X_8.
\endgather
$$
Let $\tilde\goth g$ be the affine Lie 
algebra associated with $\goth g$ spanned by the
elements $x(n)$, $x\in \goth g$, $n\in \Bbb Z$, the canonical central element
$c$ and
a derivation $d$. Set
$$\align
&\bar{B} = \{b(n) \mid  b \in B,  n \in \Bbb Z\},\\
&\bar{B}_- =\bar{B}_{<0} \cup \{X_i(0)\mid i=6, 7, 8\},\\
&\bar{B}_{<0} = \{b(n) \mid b \in B ,  n \in \Bbb Z_{<0}\},
\endalign$$
so that $\bar{B}, \bar{B}_-, \bar{B}_{<0}$ parametrize bases of the Lie
algebras 
$\bar \goth g =\hat{\goth g}/\Bbb F c$, $\tilde{\goth n}_-$ and $\tilde{\goth
g}_{<0}$
 (respectively).
We choose the order $\preccurlyeq$ on $\bar B$ defined by 
$$
b_1 (j_1)\prec b_2 (j_2) \quad \text{iff} \quad j_1 < j_2 \quad \text{or}\quad
j_1 = j_2, \quad b_1 \prec b_2. 
$$
The set of colored partitions $\Cal P (\bar B)$ is defined as the set of all
maps
$\pi : \bar B \rightarrow \Bbb N$, where $\pi(a)$ equals zero for all but
finitely
many $a \in \bar B$. Clearly $\pi$ is determined by its values $(\pi(a)
\mid a \in
\bar B)$ and we shall write $\pi$ as the monomial
$
\pi = \prod_{a\in \bar{B}} a^{\pi(a)}. 
$
We may also think of the colored partition  $\pi = \prod b_i (j_i)$ as 
$$
\pi =(b_1 (j_1),\dots, b_s (j_s)), \quad
 b_1 (j_1)\preccurlyeq\dots\preccurlyeq b_s (j_s), 
$$
where $b_i (j_i) \in \bar B$ are called the parts of $\pi$,
$|\pi |=j_1 +\dots + j_s$ the degree of $\pi$ and $\ell (\pi)=s \ge 0$ the
length of $\pi$. Since the basis elements $X_i$ are weight vectors for the
adjoint action of $\goth g$, we define the $\goth h$-weight $\wt (\pi)$ as
the sum of
the $\goth h$-weights of the $b_i$'s.
We shall think of $\pi$ as the ``plain'' partition $\sh \pi$
of the form $j_1\leq j_2\leq \dots \leq j_s$ ``colored'' with ``colors''
$b_1, b_2, \dots ,b_s$ {}from the set of colors $B=\{X_1,\dots ,X_8\}$. 
Sometimes we shall shortly say that $X_{i_1}(j_1)X_{i_2}(j_2)$ is  colored 
with colors $i_1i_2$.

Let $\pi,\pi' \in \Cal P (\bar B)$,
 $\pi=(a_1,\dots, a_s)$, $\pi'=(a' _1,\dots,a' _{s'})$. 
We extend the order $\preccurlyeq$ on $\bar B \subset \Cal P (\bar B)$ 
to the order on $\Cal P (\bar B)$ defined by 
$$
\pi \prec \pi'
$$
if $\pi \neq \pi'$ and one of the following statements holds:
\roster 
\item "{(i)}" $\ell(\pi) > \ell (\pi')$,
\item "{(ii)}" $\ell (\pi) = \ell (\pi')$,
  $\vert \pi \vert < \vert \pi' \vert$,
\item "{(iii)}"  $\ell (\pi) = \ell (\pi')$, 
 $\vert \pi \vert = \vert \pi' \vert$
and there is $i$, $\ell(\pi) \ge i \ge 1$, such that 
$\vert a_j \vert=\vert a' _j\vert$ for $\ell(\pi) \ge j > i$ and
$\vert a_i \vert < \vert a' _i \vert$,
\item"{(iv)}" $\sh \pi=\sh \pi'$ and there is $i$, $\ell(\pi) \ge i \ge 1$,
such that
$a_j = a' _j$ for $\ell(\pi) \ge j > i$ and $a_i \prec a' _i$.
\endroster

We also order the plain partitions  by the requirements (i)--(iii).
So if (i), (ii) or (iii) holds, we have that $\sh \pi \prec \sh \pi'$ and 
(hence)  $\pi \prec \pi'$.

We denote by $U_1 (\tilde \goth g)$ the quotient of the universal enveloping
algebra 
of $\tilde\goth g$ by the ideal generated by $c-1$, and we denote by
$\overline{U_1 (\tilde \goth g)}$ the completed
enveloping algebra, where  $x_i \rightarrow x$ if for each
$\tilde \goth g$-module 
$V$ of level 1 in the category $\Cal O$ and each vector $v$ in $V$ there is
$i_0$
such that $i \ge i_0$ implies $x_i v = x_{i_0} v = x v$ \cite{MP1},
\cite{MP2, 6.4}.
For $\pi \in \Cal P (\bar B)$ set
$$\align
&u(\pi)=b_1 (j_1)\dots b_s (j_s) \in U_1 (\tilde{\goth g}),\\
&U_{[\pi]} =\overline{\lspan\{u(\pi')\mid \pi' \succcurlyeq \pi\}},\\
&U_{(\pi)} =\overline{\lspan\{u(\pi')\mid \pi' \succ \pi\}},
\endalign$$
the closure taken in $\overline{U_1(\tilde \goth g)}$.
For $u \in U_{[\pi]}$, $u \notin U_{(\pi)}$, we say that $\pi$ is the
leading term
of $u$ and we write $\lt(u) = \pi$.

Recall that the generalized Verma $\tilde\goth g$-module
$N({\Lambda}_0)\cong U(\tilde\goth g_{<0})$ has the structure of vertex
operator
algebra (cf. \cite{MP2, Section 3}).
Set $r_{2\theta}=X_1(-1)X_1(-1)\bold 1$ and $R = U(\goth g)r_{2\theta}$. Then
$R$ is
an irreducible $\goth g$-module with highest weight vector $r_{2\theta}$. Set
$$
\bar R =\coprod_{n\in\Bbb Z}R(n),\quad R(n)= \{r(n) \mid r\in R \},
$$
where $Y(r,z)=\sum r_n z^{-n-1}$ denotes the vertex operator associated with
the
vector $r$ and $r(n)=r_{n+1}$ for $r\in R$. 
We think of $\bar R$ as a subset of the completed
enveloping algebra $\overline{U_1 (\tilde \goth g)}$.
Each nonzero element in $\bar R$ has a leading term. Moreover,
since $\bar R$ is generated by the adjoint action of $\goth g$ on the
elements of the form
$$
   r_{2\theta}(n)=(r_{2\theta})_{n+1}=
\sum_{i+j=n}X_1(i)X_1(j), \qquad n\in \Bbb Z,\tag 1
$$
we can describe the set $\lt (\bar R) = \{\lt(r)\mid r\in
\bar{R}\setminus\{0\}\}$ explicitly:

\proclaim{Lemma 1}
The set $\lt (\bar R)$ consists of the elements of the form
$$\gather
X_{i_1}(j)X_{i_2}(j), \text{ with colors \ } i_1i_2:
 11,21,22,31,32,33,42,43,44,51,52,53,54,55,\\
\qquad \qquad\qquad\qquad\qquad \qquad\qquad\qquad
62,64,65,66,73,74,75,76,77,85,86,87,88,\\
X_{i_1}(j-1)X_{i_2}(j), \text{ with colors \ } i_1i_2:
  11,12,13,14,15,16,17,18,22,24,26,27,28,\\
\qquad \qquad\qquad\qquad\qquad \qquad\qquad\quad
33,35,36,37,38,47,48,56,58,66,68,77,78,88,
\endgather$$
where $j\in \Bbb Z$.
\endproclaim

\demo{Proof}
For $n=2j-1$ the element 
$r_{2\theta}(2j-1)=2\cdot X_1(j\!-\!1)X_1(j)+\dots $ has the
leading term $X_1(j\!-\!1)X_1(j)$. The
adjoint action of $f_1$ on $r_{2\theta}(2j-1)$ gives the element 
$2X_3(j\!-\!1)X_1(j)+2X_1(j\!-\!1)X_3(j)+\dots$ with the leading term 
$X_1(j\!-\!1)X_3(j)$. In this way we can construct a basis of the
$\goth g$-module
$R(n)\subset \bar R$ of dimension 27. The case when $n=2j$ is treated
similarly.\qed
\enddemo

Clearly we can fix a map 
\ $\lt (\bar{R}) \rightarrow \bar R $, $\rho \mapsto r(\rho)$,
such that $\lt (r(\rho)) = \rho$.
Moreover, we will assume that this map is such that $r(\rho)$ is 
homogeneous of $\goth h$-weight $\wt (\rho)$ and degree $\vert \rho \vert$,
and that the coefficient $c$ of ``the leading term'' $X_i(n)X_j(m)$ in
``the expansion'' of $r(X_i(n)X_j(m))=cX_i(n)X_j(m)+\dots$ is chosen to be
$c=1$.
Note that $\{r(\rho)\mid \rho\in \lt (\bar R)\}$
is a basis of $\bar R$ described explicitly by Lemma 1.

Let
$$
\bar\goth g\bar\otimes\bar R = \bigoplus_{n\in \Z}\prod_{i+j=n}\goth g(i)
\otimes R(j)
$$
and use infinite sum notation to denote the elements in
$\bar\goth g\bar\otimes\bar R$ by
e.g.
$$\sum_{i+j=n}x^{(i)}(i)\otimes r^{(j)}(j),$$
$x^{(i)}\in \goth g$,
$r^{(j)}\in R$.
Since $\bar \goth g$ and $\bar R$ are loop modules,
$\bar\goth g\bar\otimes\bar R$ becomes
a $\bar\goth g$-module in the natural way. Define a linear map
$
\Psi \: \bar\goth g\bar\otimes\bar R \rightarrow \overline{U_1
(\tilde \goth g)}
$
by linear extension of
$$
\Psi \: \sum_{i+j=n}x^{(i)}(i)\otimes r^{(j)}(j) \mapsto
\sum_{i<0, i+j=n}x^{(i)}(i) r^{(j)}(j) +\sum_{i\geq 0, i+j=n}r^{(j)}(j)
x^{(i)}(i).
$$
Then $\Psi $ is not a $\bar\goth g$-module map, but it is a $\goth g$-module
map.

\proclaim{Lemma 2}
The following elements in $\bar\goth g\bar\otimes\bar R$ are for all
$n\in \Bbb Z$ 
highest weight 
vectors of $\goth g$-submodules of dimensions 64, 35, 35 and 27 (respectively):
$$\align
&q_{64}(n)=\sum_{j\in \Bbb Z} (3j-n)X_1(j)\otimes r_{2\theta}(n-j),\\
&q_{35}(n)=\bigl (f_1(-1)h_1(0)-f_1(0)h_1(-1)\bigr )\cdot q_{64}(n+1),\\
&q_{\underline {35}}(n)=\bigl (f_2(-1)h_2(0)-f_2(0)h_2(-1)\bigr )\cdot
q_{64}(n+1),\\
&q_{27}(n)=\bigl (f_2(-1)h_2(0)-f_2(0)h_2(-1)\bigr )\cdot q_{35}(n+1).
\endalign $$ 
\endproclaim
\demo{Proof}
Note first that with
$$
q'=\bigl (f_1(-1)h_1(0)-f_1(0)h_1(-1)\bigr )\cdot q,
$$ 
(a) $e_1 q=e_1(-1) q=0$ implies $e_1 q'=0$ and 
$e_2 q=e_2(-1) q=0$ implies $e_2 q'=0$, and 
(b) $e_1(-1) q=e_1(-2) q=\bigl(h_1(-2)h_1(0)-h_1(-1)h_1(-1)\bigr) q=0$ 
implies $e_1(-1) q'=0$ and
$e_2(-1) q=e_2(-2) q=0$ implies $e_2(-1) q'=0$.

It is clear that $q_{64}(n)$ is a highest weight vector. Moreover, for $i=1,2$ 
and $j\in \Bbb Z$ we have
$e_i(j) q_{64}(n)=0$ and $\bigl(h_i(-2)h_i(0)-h_i(-1)h_i(-1)\bigr)
q_{64}(n)=0$.
Now by applying (a) and (b) we see that the listed vectors are
highest weight vectors (of course, provided they are not zero, see 
Lemma 5). Since the generators $f_1$ and $f_2$ act nilpotently,
these vectors generate finite dimensional $\goth g$-modules with
dimensions given by the Weyl formula. \qed 
\enddemo

\proclaim{Proposition 3}
For all $n\in \Bbb Z$ we have
$$\align
&\Psi (q_{64}(n))=\Psi (q_{35}(n))=\Psi (q_{\underline{35}}(n))=0,\\
&\Psi (q_{27}(n))=c(n) r_{2\theta}(n)\quad \text{for some }c(n)\in \Bbb F.
\endalign$$
\endproclaim
\demo{Proof}
Since we have 
$L_{-1}X_1(-1)^3\bold 1=3X_1(-2)X_1(-1)^2\bold 1$,
where  $L_{-1}$ is a Virasoro algebra element, we can write the vertex
operator associated with this vector in two different ways:
$$
\frac d{dz}\biggl(Y(X_1(-1)\bold 1,z)Y(r_{2\theta},z)\biggr)
=3\biggl(\frac d{dz}Y(X_1(-1)\bold 1,z)\biggr) Y(r_{2\theta},z).
$$
The coefficients of this relation 
give $\Psi (q_{64}(n))=0$ for all $n\in \Bbb Z$.

For $i=1,2$ set $A_i=f_i(-1)h_i(0)-f_i(0)h_i(-1)$. Since 
$x(i)\otimes r(j) \mapsto x(i)r(j)$ and
$x(i)\otimes r(j) \mapsto r(j)x(i)$ are $\bar\goth g$-module maps, it is
clear that
$$\align
 &\Psi\biggl(A_1 \sum_{i<0,i+j=n}x^{(i)}(i)\otimes r^{(j)}(j)\biggr) = 
A_1 \Psi\biggl( \sum_{i<0,i+j=n}x^{(i)}(i)\otimes r^{(j)}(j)\biggr), \\
&\Psi\biggl(A_1 \sum_{i\geq 1,i+j=n}x^{(i)}(i)\otimes r^{(j)}(j) \biggr)= 
A_1\Psi\biggl( \sum_{i\geq 1,i+j=n}x^{(i)}(i)\otimes r^{(j)}(j) \biggr).
\endalign$$
This implies that
$$\align
&\Psi (q_{35}(n-1))=\Psi (A_1q_{64}(n))
=A_1\Psi (q_{64}(n))\\
&+\Psi\biggl(A_1\sum_{j=0}^0 (3j-n)X_1(j)\otimes r_{2\theta}(n-j)\biggr)\\
&-A_1\biggl(\sum_{j=0}^0 (3j-n)r_{2\theta}(n-j)X_1(j)\biggr).
\endalign$$
Now note that the last two terms consist of products of $X_1(i)$ or
$f_1\cdot X_1(i)$
with $r_{2\theta}(i)$ or $f_1\cdot r_{2\theta}(i)$, $i\in \Bbb Z$, differing
only in 
the order they are written. Since these elements commute (cf. (1)), the last
two terms cancel. This together with $\Psi (q_{64}(n))=0$ gives
$\Psi (q_{35}(n-1))=0$. The relation $\Psi (q_{\underline{35}}(n))=0$ is
proved in
a similar way.

Arguing as before we get
$$\align
&\Psi (q_{27}(n-2))=\Psi (A_2A_1q_{64}(n))
=A_2A_1\Psi (q_{64}(n))\\
&+\Psi\biggl(A_2A_1\sum_{j=0}^1 (3j-n)X_1(j)\otimes r_{2\theta}(n-j)\biggr)\\
&-A_2A_1\biggl(\sum_{j=0}^1 (3j-n)r_{2\theta}(n-j)X_1(j)\biggr).
\endalign$$
Now note that the last two terms consist of products of
$X_1(i)$, $f_1\cdot X_1(i)$,
$f_2\cdot X_1(i)$ or $f_2f_1\cdot X_1(i)$
with $r_{2\theta}(i)$, $f_1\cdot r_{2\theta}(i)$, 
$f_2\cdot r_{2\theta}(i)$ or $f_2f_1\cdot r_{2\theta}(i)$, $i\in \Bbb Z$,
differing 
only in the order they are written. Since $\bar R$ is a loop module, that is
$[x(i), r(j)]=(x\cdot r)(i+j)$, we can commute these elements and the last
two terms cancel except for an element in $\bar R$ of $\goth h$-weight
$2\theta$
and degree $n-2$. This together with $\Psi (q_{64}(n))=0$ gives that
$\Psi (q_{27}(n-2))$ is proportional to $ r_{2\theta}(n-2)$. \qed
\enddemo

Recall that $N^1(\Lambda_0)=\bar R N(\Lambda_0)$ is the maximal
$\tilde\goth g$-submodule
of $N(\Lambda_0)$.
The elements $r(n)=r_{n+1}\in\bar R$ annihilate the basic
$\tilde\goth g$-module
$L(\Lambda_0)=N(\Lambda_0)/N^1(\Lambda_0)$, so we call them relations. 
On the other hand $\Psi (q_{64}(n))=0$ reads
$$
\sum_{j<0}(3j-n)X_1(j) r_{2\theta}(n-j) +\sum_{j\geq 0}(3j-n)r_{2\theta}(n-j)
X_1(j)
=0 \quad \text{on \ }N(\Lambda_0),
$$
so we shall sometimes say that  $\Psi (q_{64}(n))=0$, 
or $q_{64}(n)\in\bar\goth g\bar\otimes\bar R$, are
relations among relations. If we extend 
$\Psi \: \bar\goth g\bar\otimes\bar R +\Bbb F \otimes \bar R
\rightarrow \overline{U_1 (\tilde \goth g)}$ by 
$\Psi \bigl(1\otimes r(n)\bigr)=r(n)$, then we can write
$\Psi (q_{27}(n))=c(n) r_{2\theta}(n)$ as
 $\Psi \bigl(q_{27}(n)-c(n) \otimes r_{2\theta}(n)\bigr)=0$.
By abuse of notation we shall sometimes
write or think $q_{27}\in \text{ker}\,\Psi$.

For colored partitions $\kappa$, $\rho$ and $\pi=\kappa\rho$ we shall write
$\kappa=\pi/\rho$ and $\rho\subset\pi$. We shall say that $\rho\subset\pi$ is
an
embedding (of $\rho$ in $\pi$).
For an embedding $\rho\subset\pi$, where $\rho\in\lt (\bar R)$, we define the
element $u(\rho\subset\pi)$ in $ \overline{U_1(\tilde \goth g)}$ by
$$
u(\rho\subset\pi)=\cases u(\pi/\rho)r(\rho) &\text{if \ } |\rho
|>|\pi/\rho |,\\
r(\rho)u(\pi/\rho) &\text{if \ } |\rho |\leq |\pi/\rho |.\endcases \tag 2
$$
It is clear that each nonzero $q\in\bar\goth g\bar\otimes\bar R$ 
(or $q\in\bar\goth g\bar\otimes\bar R +\Bbb F\otimes \bar R$) can be written
in the form
$$
q =  \sum_{\rho' \subset \pi} C_{\rho',\pi}\, u(\pi/\rho') \otimes
r (\rho') 
 + \sum \Sb \rho' \subset \pi'\\ \pi' \succ \pi \endSb C_{\rho', \pi'}\,
 u(\pi'/\rho') \otimes r (\rho')  
$$
for some $\pi$ of length 3 (or $\leq 3$) and some 
$C_{\rho',\pi}, C_{\rho',\pi'} \in \Bbb F$, 
where at least one coefficient $C_{\rho',\pi}$
is nonzero. We shall say that $\pi$ is the leading term of $q$ and we write
$\pi=\lt (q)$.
The assumption $q \in \ker \Psi$, i.e.
$$
 \sum_{\rho' \subset \pi} C_{\rho',\pi} u(\rho' \subset \pi) +
 \sum \Sb \rho' \subset \pi'\\ \pi' \succ \pi \endSb D_{\rho', \pi'}
 u(\rho' \subset \pi')  = 0,
$$
implies that $\sum_{\rho' \subset \pi} C_{\rho',\pi} = 0$
(cf. [MP2, Lemma 6.4.1]). Hence for such $q$ there
must be at least two (different) embeddings $\rho' \subset \pi$, 
$\rho'' \subset \pi$ with
the corresponding coefficients being nonzero.
For a colored partition $\pi$ of length 3 set 
$$
N(\pi)=\max \{\#\Cal E(\pi)\!-\!1,\,0\},\quad 
\Cal E(\pi)=\{\rho\in \lt (\bar R) \mid \rho\subset\pi\}.
$$
It is clear that $0\leq N(\lt (q))\leq 2$ for $q\in\bar\goth g\bar\otimes\bar
R$.
The argument above shows that $N(\lt (q))\geq 1$ for $q\in ker \Psi$.

\proclaim{Lemma 4}
Let $Q\subset ker \Psi$ be a subspace of dimension $n$. Then
$$
\sum_{\pi\in\lt (Q)} N(\pi)\geq n,
$$
where $\lt (Q)=\{\lt (q)\mid q\in Q\backslash \{0\}\}$.
\endproclaim
\demo{Proof}
Let $Q_{[\pi]}=\{q\in Q \mid \lt (q) \succcurlyeq \pi\}$ and
$Q_{(\pi)}=\{q\in Q \mid \lt (q) \succ \pi\}$. Let 
$\dim Q_{[\pi]}/Q_{(\pi)}=m(\pi)=m\leq n$,
$m\geq 1$ and let $\rho_1\subset\pi,\dots ,\rho_s\subset\pi$
(where $s=s(\pi)$) 
be all possible
embeddings in $\pi$. Let $\pi^*$ be such that $Q_{(\pi)}=Q_{[\pi^*]}$.
Then we can write a basis of $Q_{[\pi]}/Q_{(\pi)}$ in the form
$$\aligned
c_{11}u(\pi/\rho_1)\otimes r(\rho_1)+&\dots+c_{1s}u(\pi/\rho_s)\otimes
r(\rho_s)
+v_1+Q_{(\pi)},\\
c_{21}u(\pi/\rho_1)\otimes r(\rho_1)+&\dots+c_{2s}u(\pi/\rho_s)\otimes
r(\rho_s)
+v_2+Q_{(\pi)},\\
&\dots\\
c_{m1}u(\pi/\rho_1)\otimes r(\rho_1)+&\dots+c_{ms}u(\pi/\rho_s)\otimes
r(\rho_s)
+v_m+Q_{(\pi)},
\endaligned $$
where the vectors $v_i$ are of the form 
$$
v_i=\sum_{\pi\prec \pi'\prec \pi^*}\sum_{\rho\in\Cal E(\pi')}
d_{i,\rho,\pi'}\,u(\pi'/\rho)\otimes r(\rho)+
\sum_{\pi^*\preccurlyeq \pi'}\sum_{\rho\in\Cal E(\pi')}
e_{i,\rho,\pi'}\,u(\pi'/\rho)\otimes r(\rho).
$$
Assume that $\text{rank\,}(c_{ij})<m$. Then the rows are linearly dependent.
By taking a nontrivial linear combination of the basis elements 
we get a vector in $Q$ of the form
$$
v=\sum_{\pi\prec \pi'\prec \pi^*}\sum_{\rho\in\Cal E(\pi')}
d_{\rho,\pi'}\,u(\pi'/\rho)\otimes r(\rho)+
\sum_{\pi^*\preccurlyeq \pi'}\sum_{\rho\in\Cal E(\pi')}
e_{\rho,\pi'}\,u(\pi'/\rho)\otimes r(\rho).
$$
The coefficients $d_{\rho,\pi'}$, $\rho\in\Cal E(\pi')$, $\pi\prec \pi'\prec
\pi^*$,
must be zero since otherwise $Q_{(\pi)}\neq Q_{[\pi^*]}$. But then
$v\in Q_{(\pi)}=Q_{[\pi^*]}$ and our nontrivial linear combination of
basis elements is zero in $Q_{[\pi]}/Q_{(\pi)}$, a contradiction.

Hence the rank of the matrix $(c_{ij})$ is $m$ and we have $s\geq m$.
Assume that
$s=m$. Then the matrix $(c_{ij})$ is regular and the set of
vectors of the form $u(\pi/\rho_1)\otimes r(\rho_1)+v_1+Q_{(\pi)},\dots ,
u(\pi/\rho_s)\otimes r(\rho_s)+v_m+Q_{(\pi)}$ (vectors $v_i$ as above) 
is a basis of $Q_{[\pi]}/Q_{(\pi)}$.
In particular, we have a vector $u\in Q$ of the form 
$u= u(\pi/\rho_1)\otimes r(\rho_1)+v_1$ such that $\Psi (u)\neq 0$, a
contradiction. 

Hence $s>m$, i.e. $N(\pi)= s(\pi)-1\geq m(\pi)$. Since 
$\dim Q=\sum m(\pi)$, the lemma follows.
\qed
\enddemo
\remark{Remark}
Let us keep the notation {}from the proof of Lemma 4 and let us
assume that $N(\pi)=m(\pi)$. Then there are altogether $s=m+1$ possible
embeddings $\rho_1\subset\pi,\dots ,\rho_s\subset\pi$ and the rank of
$(c_{ij})$
is $m=s-1$. Let us assume that the first $s-1$ columns of $(c_{ij})$ are
linearly
independent. Then for each $1\leq i\leq s-1$ there is a vector in $Q$ of
the form
$$
u(\pi/\rho_i)\otimes r(\rho_i)+d_iu(\pi/\rho_s)\otimes r(\rho_s)
+\sum_{\pi\prec \pi'}\sum_{\rho\in\Cal E(\pi')}
d_{i,\rho,\pi'}\,u(\pi'/\rho)\otimes r(\rho)
$$
for some $d_i, d_{i,\rho,\pi'}\in \Bbb F$. As was shown above, the assumption
$Q\subset ker \Psi$ implies $d_i\neq 0$. But then for
each $i,j\in\{1,\dots ,s\}$
there are $d_{ij}, d_{i,j,\rho,\pi'}\in \Bbb F$ such that
$$
u(\pi/\rho_i)\otimes r(\rho_i)+d_{ij}u(\pi/\rho_j)\otimes r(\rho_j)
+\sum_{\pi\prec \pi'}\sum_{\rho\in\Cal E(\pi')}
d_{i,j,\rho,\pi'}\,u(\pi'/\rho)\otimes r(\rho) \in ker \Psi . \tag 3
$$
\endremark
As above we denote with by dot the adjoint action of $U(\goth g)$ and set 
$$
Q_A(n)=U(\goth g) \cdot q_{A}(n) \text{ \ for \ } n\in \Bbb Z 
\text{ \ and \ } A=27,35,\underline{35},64.
$$

\proclaim{Lemma 5}
For $j\in \Bbb Z$ we have:
$$\align
&\sh (\lt (q))=j^3\quad \text{for \ }q\in  Q_{27}(3j)\oplus Q_{35}(3j)\oplus 
Q_{\underline{35}}(3j),\\
&\sh (\lt (q))=(j\!-\!1)j(j\!+\!1)\quad \text{for \ }q\in  Q_{64}(3j),\\ 
&\sh (\lt (q))=(j\!-\!1)j^2 \text{\ for }q\in Q_{27}(3j\!-\!1)\oplus 
Q_{35}(3j\!-\!1)\oplus Q_{\underline{35}}(3j\!-\!1)\oplus Q_{64}(3j\!-\!1),\\
&\sh (\lt (q))=(j\!-\!1)^2j\text{\ for }q\in Q_{27}(3j\!-\!2)\oplus 
Q_{35}(3j\!-\!2)\oplus Q_{\underline{35}}(3j\!-\!2)\oplus Q_{64}(3j\!-\!2). 
\endalign$$
\endproclaim
\demo{Proof}
Note first that each $q_A(n)$, $A=27,35,\underline{35},64$, is a sum of
elements
of the form $C_{\rho', \pi'}\, u(\pi'/\rho') \otimes r (\rho')$ with 
$|\pi'|=n$ and $\ell (\pi')=3$
(and additional terms with $\ell (\pi')=2$ appearing in $q_{27}$), and that
for
such colored partitions $\sh (\pi')=abc$, where $a\leq b\leq c$,
$a+b+c=|\pi'|$. The smallest possible shapes are
$$
\alignat2
&j^3\prec (j-1)j(j+1)\prec \dots \qquad&&\text{when \ } |\pi'|=3j,\\
& (j-1)j^2\prec \dots \qquad&&\text{when \ } |\pi'|=3j-1,\\
& (j-1)^2j\prec \dots \qquad&&\text{when \ } |\pi'|=3j-2.
\endalignat$$
In the case of $q_{64}(3j)$ we see that for $\sh (\pi')=j^3$ the coefficient
$C_{\rho', \pi'}$ equals $3j-3j=0$, and that for $\sh (\pi')=(j-1)j(j+1)$ the
corresponding coefficients are not zero. Hence 
$\sh \bigl(\lt (q_{64})\bigr)=(j-1)j(j+1)$. Since the projection
$$\align
&\phi \: Q_{64}(3j) \rightarrow \goth g(j-1)\otimes R(2j+1) 
+ \goth g(j+1)\otimes R(2j-1),\\
&\phi \: U(\goth g)q_{64}(3j)\rightarrow U(\goth g)
\bigl(-3X_1(j-1)\otimes r_{2\theta}(2j+1)+
3X_1(j+1)\otimes r_{2\theta}(2j-1)\bigl)
\endalign$$
is a $\goth g$-module map, for each $q\in U(\goth g)q_{64}(3j)$, $q\neq 0$,
we have
$\phi (q)\neq 0$. Hence  $\sh \bigl(\lt (q)\bigr)=(j-1)j(j+1)$.

In the case of $q_{64}(3j-1)$ we see that for $\sh (\pi')=(j-1)j^2$ the
coefficients
$C_{\rho', \pi'}$ equal $3j-3-3j+1=-2\neq 0$ and $3j-3j+1=1\neq 0$, 
and hence 
$\sh \bigl(\lt (q_{64})\bigr)=(j-1)j^2$. Since the projection
$$
\phi \: U(\goth g)q_{64}(3j-1)\rightarrow U(\goth g)
\bigl(-2X_1(j-1)\otimes r_{2\theta}(2j)+
X_1(j)\otimes r_{2\theta}(2j-1)\bigl)
$$
is a $\goth g$-module map, for each $q\in U(\goth g)q_{64}(3j-1)$, $q\neq 0$,
we have
$\phi (q)\neq 0$. Hence  $\sh \bigl(\lt (q)\bigr)=(j-1)j^2$.

The other cases  $q_A(n)$, $A=27,35,\underline{35}$, are similar (except for
more
complicated expressions for $q_A$): we have to show that 
at least one coefficient $C_{\rho', \pi'}\neq 0$ for each of the shapes
$j^3$, $(j-1)j^2$ and $(j-1)^2j$.\qed

\enddemo

Note that Lemma 4 and Lemma 5 imply that:
\roster
\item"{}" $\sum_{\sh (\pi)=j^3}N(\pi)\geq 27+ 35+35=97$, 
\item"{}" $\sum_{\sh (\pi)=(j\!-\!1)j(j\!+\!1)}N(\pi)\geq 64$,
\item"{}"  $\sum_{\sh (\pi)=(j\!-\!1)j^2 }N(\pi)\geq 27+35+35+64=161$,
\item"{}"  $\sum_{\sh (\pi)=(j\!-\!1)^2j}N(\pi)\geq   27+35+35+64=161$.
\endroster
By direct counting we see the following:
\proclaim{Lemma 6}
\roster
\item"{}" $\sum_{\sh (\pi)=j^3}N(\pi)=97$, 
\item"{}" $\sum_{\sh (\pi)=(j\!-\!1)j(j\!+\!1)}N(\pi) =64$,
\item"{}"  $\sum_{\sh (\pi)=(j\!-\!1)j^2 }N(\pi)=162$,
\item"{}"  $\sum_{\sh (\pi)=(j\!-\!1)^2j}N(\pi)=162$.
\endroster
\endproclaim
Lemmas 4, 5 and 6 imply that for all but two $\goth h$-weight spaces
$Q(n)_\mu$ we have
$$
\sum_{ \ell (\pi)=3, |\pi|=n, \wt (\pi)=\mu}N(\pi)=\dim Q(n)_\mu, \tag 4
$$
where $n=3j,3j\!-\!1,3j\!-\!2$ and 
$Q(n)=Q_{27}(n)\oplus Q_{35}(n)\oplus Q_{\underline{35}}(n)
\oplus Q_{64}(n)$. The next lemma
describes precisely these two exceptions. There we shall use the notation
$$\spreadlines{-1.15ex}\alignat 2
&\sq\sq \quad &&\circ\bullet   3\quad    3\quad \circ    5\\
&\sq \quad    &&\bullet\bullet 5\quad    4\quad \bullet  3\\
&\sq \quad    &&\bullet\circ   1\quad    1\quad \bullet  1,
\endalignat
$$
where the Young diagram represents the plain partition $(-2)(-1)^2$, we shall
think of $(j-1)j^2$ in general, the numbers 351, 341 and 531 represent
colorings.
So
we have listed three colored partitions ({}from left to right):
$X_3(j-1)X_5(j)X_1(j)\prec  X_3(j-1)X_4(j)X_1(j)\prec X_5(j-1)X_3(j)X_1(j)$.
Moreover, with circs and bullets we denote all possible embeddings
$X_5(j)X_1(j)\subset X_3(j-1)X_5(j)X_1(j)$,
$X_3(j-1)X_5(j)\subset X_3(j-1)X_5(j)X_1(j)$ and
$X_3(j)X_1(j)\subset X_5(j-1)X_3(j)X_1(j)$ for the listed colored
partitions (cf. Lemma 1).

\proclaim{Lemma 7}
(a) For $\mu=\alpha_1+2\alpha_2$ and $j\in \Bbb Z$ there are the following 10 
partitions $\pi$ such that $\sh (\pi)=(j-1)j^2$ and $\wt (\pi)=\mu$:
$$\spreadlines{-1.15ex}\alignat 2
&\sq\sq \quad &&\bullet\bullet\circ 1\quad \bullet\bullet\circ 1\quad \circ
2\quad 
\circ\bullet 3 \quad \bullet\bullet 1\quad \circ\bullet 3\quad 
 3\quad \circ 5\quad \circ 4\quad \circ 7 \\
&\sq \quad    &&\circ\bullet\bullet 5\quad \circ\bullet\bullet 4\quad
\bullet 3\quad 
\bullet\bullet 3\quad \circ\bullet 7\quad \bullet\bullet 5\quad   
 4\quad \bullet 3\quad \bullet 3\quad \bullet 1 \\
&\sq \quad    &&\bullet\circ\bullet 3\quad \bullet\circ\bullet 3\quad
\bullet 3\quad  
\bullet\circ 2\quad \bullet\circ 1\quad  \bullet\circ 1\quad 
 1\quad \bullet 1\quad \bullet 1\quad \bullet 1.
\endalignat
$$
Moreover, for such partitions $\sum N(\pi)=7=\dim Q(3j-1)_\mu +1$.

(b) For $\mu=-\alpha_1-2\alpha_2$ and $j\in \Bbb Z$ there are the following 10 
partitions $\pi$ such that $\sh (\pi)=(j-1)^2j$ and $\wt (\pi)=\mu$:
$$\spreadlines{-1.15ex}\alignat 2
&\sq\sq \quad &&\bullet\bullet\circ 6\quad \bullet\bullet\circ 6\quad
\circ\bullet  8\quad
\bullet 6 \quad \circ\bullet 7\quad \circ\bullet 8\quad 
 8\quad \bullet 8\quad \bullet 8\quad \bullet 8 \\
&\sq\sq \quad &&\circ\bullet\bullet 5\quad \circ\bullet\bullet 4\quad
\bullet\circ 2\quad
\bullet 6\quad \bullet\bullet 6\quad \bullet\bullet 5\quad   
 4\quad \bullet 6\quad \bullet 6\quad \bullet 8 \\
&\sq \quad    &&\bullet\circ\bullet 8\quad \bullet\circ\bullet 8\quad
\bullet\bullet 8\quad
\circ 7\quad \bullet\circ 6\quad  \bullet\circ 6\quad 
 6\quad \circ 5\quad \circ 4\quad \circ 2.
\endalignat
$$
Moreover, for such partitions $\sum N(\pi)=7=\dim Q(3j-2)_\mu +1$.
\endproclaim

 Let $\rho_1 \subset \pi$, $\rho_2 \subset \pi$ be two embeddings,
$\rho_1, \rho_2 \in
\lt(\bar{R})$. We would like to construct a relation among relations of
the form
$$
u(\rho_1 \subset \pi) \in \Bbb F^\times u(\rho_2 \subset \pi) + 
\overline {\lspan\{u(\rho \subset \pi')\mid
 \rho \in \lt (\bar{R}), \rho \subset \pi', \pi \prec \pi'\}}
\tag{5}
$$
(the closure taken in $\overline{U_1(\tilde{\goth g})}$).
If the colored partition $\pi$ is such that
$\ell (\pi)=3$ and that for $\mu=\wt (\pi)$ the relation (4)
holds, then the proof of Lemma 4 shows that $s(\pi)=m(\pi)+1$ and that there is
an element $q\in \text{ker}\,\Psi$ of the form (see (3))
$$
q=u(\pi/\rho_1)\otimes r(\rho_1)+c\, u(\pi/\rho_2)\otimes r(\rho_2)+
\sum_{\pi\prec \pi'}\sum_{\rho\in\Cal E(\pi')}
c_{\rho,\pi'}\,u(\pi'/\rho)\otimes r(\rho),
$$
and this implies (5). However, the equality (4) does not hold for all
weights $\mu$.
In the next two lemmas we identify embeddings which do not appear in relations 
among relations of the form (5):

\proclaim{Lemma 8} Let $j\in\Bbb Z$. Then
$$\gather
\lt \bigl (X_3(j\!-\!1)r(X_5(j)X_1(j))\!-\!r(X_3(j\!-\!1)X_5(j))X_1(j)\bigr )=
X_3(j\!-\!1)X_4(j)X_1(j),\\
\lt \bigl (r(X_8(j)X_5(j))X_6(j\!+\!1)\!-\!X_8(j)r(X_5(j)X_6(j\!+\!1))\bigr )=
X_8(j)X_4(j)X_6(j\!+\!1).
\endgather$$
\endproclaim

\demo{Proof}
By calculating the first few terms of $r(X_5(j)X_1(j))$ 
and $r(X_3(j\!-\!1)X_5(j))$ we get
$$\align
&X_3(j\!-\!1)r(X_5(j)X_1(j))\!-\!r(X_3(j\!-\!1)X_5(j))X_1(j)\\
&=X_3(j\!-\!1)X_5(j)X_1(j)+X_3(j\!-\!1)X_4(j)X_1(j))+\dots\\
&-\bigl (X_3(j\!-\!1)X_5(j)X_1(j)+X_5(j\!-\!1)X_3(j)X_1(j)+\dots\bigr )\\
&=X_3(j\!-\!1)X_4(j)X_1(j))-X_5(j\!-\!1)X_3(j)X_1(j)+\dots ,
\endalign$$
where the dots represent sums of terms of higher shape. The other equality
is proved similarly. \qed
\enddemo

\proclaim{Lemma 9}
Let $j\in \Bbb Z$. Then there is no element $q$ in $\text{ker}\,\Psi$ such that
the leading term $\lt (q)$ is either $X_3(j-1)X_5(j)X_1(j)$ or
$X_8(j-1)X_5(j-1)X_6(j)$.
\endproclaim

\demo{Proof}
Set $\pi=X_3(j-1)X_5(j)X_1(j)$ and let $q\in \text{ker}\,\Psi$ be of the form
$$\align
q=&X_3(j\!-\!1)\otimes r(X_5(j)X_1(j))\!-\!X_1(j)\otimes
r(X_3(j\!-\!1)X_5(j))\\
&+\sum_{\pi\prec \pi'}\sum_{\rho\in\Cal E(\pi')}
c_{\rho,\pi'}\,u(\pi'/\rho)\otimes r(\rho).
\endalign$$ {}From Lemma 7(a) we see that 
$X_3(j-1)X_5(j)X_1(j)\prec  X_3(j-1)X_4(j)X_1(j)\prec X_5(j-1)X_3(j)X_1(j)$
and that $\Cal E(X_3(j-1)X_4(j)X_1(j))=\emptyset$. Hence
$$\align
\Psi (q)=&X_3(j\!-\!1) r(X_5(j)X_1(j))\!-\!X_1(j) r(X_3(j\!-\!1)X_5(j))\\
&+\sum_{X_3(j-1)X_4(j)X_1(j)\prec \pi'}\sum_{\rho\in\Cal E(\pi')}
c_{\rho,\pi'}\,u(\pi'/\rho) r(\rho).
\endalign$$
Now Lemma 8 implies that the leading term
$\lt (\Psi (q))=X_3(j-1)X_4(j)X_1(j)$,
and in particular this implies that $\Psi (q)\neq 0$, a contradiction.

The case of  $X_8(j-1)X_5(j-1)X_6(j)$ is proved similarly. \qed
\enddemo

We may conclude our discussion with the following:
\proclaim{Proposition 10}
Let $\pi$ be a colored partition of length $3$ and let $\rho_1 \subset \pi$, 
$\rho_2 \subset \pi$ be two embeddings, $\rho_1, \rho_2 \in\lt(\bar{R})$. Let
$$
\pi\neq X_3(j-1)X_5(j)X_1(j), \quad \pi\neq X_8(j-1)X_5(j-1)X_6(j).
$$
Then the relation (5) holds:
$$
u(\rho_1 \subset \pi) \in \Bbb F^\times u(\rho_2 \subset \pi) + 
\overline {\lspan\{u(\rho \subset \pi')\mid
 \rho \in \lt (\bar{R}), \rho \subset \pi', \pi \prec \pi'\}}.
$$
\endproclaim
\demo{Proof}
We have seen that the relation (5) holds for all cases except when
$\wt (\pi)=\alpha_1+2\alpha_2$ or $\wt (\pi)=-\alpha_1-2\alpha_2$. So
consider the first case: all possible $\pi$'s are listed in Lemma 7(a),
altogether five of them allow two or more than two embeddings, say
$$
\pi_1\prec\pi_2\prec\pi_3\prec\pi_4\prec\pi_5 = X_3(j-1)X_5(j)X_1(j),
$$
and 
$$
\sum_{i=1}^4 N(\pi_i) = 6.
$$
Set $Q=Q(3j-1)_{\alpha_1+2\alpha_2}$.
Since $Q\subset\ker\Psi$, Lemma 9 implies that
$Q_{[\pi_5]}=0$. Since 
$Q=Q_{[\pi_1]}\supset \dots \supset Q_{[\pi_5]}=Q_{(\pi_4)}=0$, we have
$$
\dim Q=\sum_{i=1}^4 \dim\left (Q_{[\pi_i]}/Q_{(\pi_i)} \right )=
\sum_{i=1}^4 m(\pi_i) =6.
$$ {}From the proof of Lemma 4 we see that $N(\pi_i)= s(\pi_i)-1\geq m(\pi_i)$.
Hence $\sum N(\pi_i)= \sum m(\pi_i)$ implies that $N(\pi_i)=m(\pi_i)$
for $i=1, 2, 3, 4$, which,
as remarked, implies the relation (5) for $\pi=\pi_1,\dots ,\pi_4$.

The case when $\wt (\pi)=-\alpha_1-2\alpha_2$ is proved similarly.
 \qed
\enddemo

Set
$$
\Cal R=\lt (\bar R)\cup \{ X_3(j-1)X_4(j)X_1(j), X_8(j-1)X_4(j-1)X_6(j)\mid
j\in \Bbb Z\}.
$$
We extend the map 
$\lt (\bar{R}) \rightarrow \bar R $, $\rho \mapsto r(\rho)$, to the map
$\Cal R \rightarrow \bar R U(\tilde\goth g)$ by
$$\align
&r\bigl (X_3(j\!-\!1)X_4(j)X_1(j) \bigr )=
X_3(j\!-\!1)r(X_5(j)X_1(j))\!-\!r(X_3(j\!-\!1)X_5(j))X_1(j),\\
&r\bigl (X_8(j)X_4(j)X_6(j\!+\!1)\bigr )=
r(X_8(j)X_5(j))X_6(j\!+\!1)\!-\!X_8(j)r(X_5(j)X_6(j\!+\!1)),
\endalign$$
so that $\lt (r(\rho)) = \rho$.
We still have that  $r(\rho)$ is 
homogeneous of $\goth h$-weight $\wt (\rho)$ and degree $\vert \rho \vert$,
and that the coefficient $c$ of ``the leading term'' $u(\rho)$ in
``the expansion'' of $r(\rho)=c\,u(\rho)+\dots$ 
is chosen to be $c=1$.
For an embedding $\rho\subset\pi$, where $\rho\in\Cal R$, we define the
element $u(\rho\subset\pi)$ in $ \overline{U_1(\tilde \goth g)}$ the 
same way as before (cf. (2)). For $\pi \in \Cal P(\bar B)$ set
$$
\Cal R_{(\pi)} = \overline {\lspan\{u(\rho \subset \pi')\mid
 \rho \in \Cal R, \rho \subset \pi', \pi \prec \pi'\}},
$$
the closure taken in $\overline{U_1(\tilde{\goth g})}$.

Note that for $\pi = X_3(j-1)X_5(j)X_1(j)$, $\rho_1= X_5(j)X_1(j)$,
$\rho_2= X_3(j-1)X_5(j)$, we have 
$$
u(\rho_1 \subset \pi) \in  u(\rho_2 \subset \pi) +\Cal R_{(\pi)}
$$
since $u(\rho_1 \subset \pi) =r(X_3(j-1)X_4(j)X_1(j))+ u(\rho_2 \subset \pi) +
\text{shorter terms}$ (cf. (2)).
So our construction gives
$$
u(\rho_1 \subset \pi) \in \Bbb F^\times u(\rho_2 \subset \pi) + \Cal R_{(\pi)}
$$
for any two embeddings $\rho_1 \subset \pi$, $\rho_2 \subset \pi$ when
$\rho_1, \rho_2 \in\lt(\bar{R})$ and  $\ell (\pi)=3$.

Our main result about relations among relations is the following:
\proclaim{Theorem 11}
 Let $\rho_1 \subset \pi$, $\rho_2 \subset \pi$ be two embeddings, $\rho_1,
\rho_2 \in
\Cal R$. Then 
$$
u(\rho_1 \subset \pi) \in \Bbb F^\times u(\rho_2 \subset \pi) + \Cal R_{(\pi)}.
\tag{6}
$$

\endproclaim

We prove the theorem in several steps, basically following the ideas used
in the case $\goth g=\goth {sl}(2,\Bbb F)$:

The case when $\rho_1\rho_2\subset \pi$ is easy: We ``expand'' the
product $r(\rho_1)r(\rho_2)$ in two different ways and get (6). For example,
for $\pi=X_3(j-1)X_5(j)X_5(j)X_1(j)$ we have
$$\align
&r(\rho_1)=r\bigl(X_3(j-1)X_5(j)\bigr)= X_3(j-1)X_5(j) +X_5(j-1)X_3(j)
+\dots,\\
&r(\rho_2)=r\bigl(X_5(j)X_1(j)\bigr)= X_5(j)X_1(j) +X_4(j)X_1(j) +\dots,
\endalign$$
where the dots represent sums of terms of higher shape. So we can
``expand'' the
product $r(\rho_1)r(\rho_2)$ in two ways and get
$$\aligned
&r(\rho_1)X_5(j)X_1(j)-X_3(j-1)X_5(j)r(\rho_2)\\
&=-r(\rho_1)X_4(j)X_1(j)+X_5(j-1)X_3(j)r(\rho_2)+\dots\\
&\in\Cal R_{(X_3(j-1)X_5(j)X_5(j)X_1(j))},
\endaligned$$
where dots represent a sum of terms of higher shape. For later purposes note
that $\lt \bigl(X_5(j-1)X_3(j)r(\rho_2)\bigr)=X_5(j-1)X_5(j)X_3(j)X_1(j)
\succ X_3(j-1)X_5(j)X_4(j)X_1(j)$ and hence
$$\aligned
&r(\rho_1)X_5(j)X_1(j)-X_3(j-1)X_5(j)r(\rho_2)\\
&\in -r(\rho_1)X_4(j)X_1(j)+
\Cal R_{(X_3(j-1)X_5(j)X_4(j)X_1(j))},
\endaligned\tag 7$$

The case when $\rho_1, \rho_2 \in \lt (\bar R)$ follows {}from the discussion
above.
In particular, relations among relations (5) hold for all $\pi$ of length 3
except
$X_3(j-1)X_5(j)X_1(j)$ or $X_8(j-1)X_5(j-1)X_6(j)$.

The remaining case is when $\rho_1\in\Cal R$ ``intersects'' 
$\rho_2\in \Cal R\backslash \lt (\bar R)$:
Define $\pi \cup \nu$ and $\pi \cap \nu$ by
$$
(\pi \cup \nu)(a) = \max \{\pi(a), \nu(a)\},\quad 
(\pi \cap \nu)(a) = \min \{\pi(a), \nu(a)\},
$$
and denote the partition $1=\prod_{a\in \bar{B}} a^0$ with no parts
and length 0 as $\varnothing$. By using Lemma 1 we see the following:

\proclaim{Lemma 12}
Let $\pi$ be a colored partition, $\rho_1\in\Cal R$ and 
$\rho_2$ either $X_3(j\!-\!1)X_4(j)X_1(j)$ or $X_8(j-1)X_4(j-1)X_6(j)$,
$j\in\Bbb Z$. 
Assume that
$\rho_1, \rho_2 \subset \pi$,
$\pi=\rho_1 \cup \rho_2$, $\rho_1 \cap \rho_2 \neq \varnothing$ and
 $\ell(\pi)\ge 4$. Then $\pi$ is one of the following colored
partitions :
$$\allowdisplaybreaks \alignat2
&X_a^\bullet (j-2)X_3^\bullet (j-1)X_4(j)X_1(j), \quad &&a=3,1, \\
&X_a^\bullet (j-1)X_3^\bullet (j-1)X_4(j)X_1(j),
\quad &&a=7,5,4,3,2,1,\\                     
&X_a^\bullet (j-1)X_3(j-1)X_4^\bullet (j)X_1(j), \quad &&a=2,1, \\
&X_a^\bullet (j-1)X_3(j-1)X_4(j)X_1^\bullet(j),
\quad &&a=1, \\                                
&X_3^\bullet (j-1)X_4(j)X_1(j)X_a^\bullet (j),
\quad &&a=8,7,6,5,3, \\                               
&X_3(j-1)X_4^\bullet (j)X_1(j)X_a^\bullet (j),
\quad &&a=7,6,5,4,3,2, \\
&X_3(j-1)X_4(j)X_1^\bullet (j)X_a^\bullet (j), \quad &&a=5,3,2,1,   \\
&X_3(j-1)X_4^\bullet (j)X_1(j)X_a^\bullet (j+1),
\quad &&a=8,7, \\                                                        
&X_3(j-1)X_4(j)X_1^\bullet (j)X_a^\bullet (j+1),
\quad &&a=8,\dots,1, \\                
&X_a^\bullet (j-2)X_8^\bullet (j-1)X_4(j-1)X_6(j),
\quad &&a=8,\dots ,1,   \\  
&X_a^\bullet (j-2)X_8(j-1)X_4^\bullet (j-1)X_6(j),  \quad &&a=2,1,   \\ 
&X_a^\bullet (j-1)X_8^\bullet (j-1)X_4(j-1)X_6(j),  \quad &&a=8,7,6,5,   \\
&X_a^\bullet (j-1)X_8(j-1)X_4^\bullet (j-1)X_6(j),
\quad &&a=7,6,5,4,3,2,   \\ 
&X_a^\bullet (j-1)X_8(j-1)X_4(j-1)X_6^\bullet (j),  \quad &&a=6,5,3,2,1,   \\ 
&X_8^\bullet (j-1)X_4(j-1)X_6(j)X_a^\bullet (j),  \quad &&a=8,   \\
&X_8(j-1)X_4^\bullet (j-1)X_6(j)X_a^\bullet (j),  \quad &&a=8,7,   \\ 
&X_8(j-1)X_4(j-1)X_6^\bullet (j)X_a^\bullet (j),  \quad &&a=8,7,6,5,4,2,   \\ 
&X_8(j-1)X_4(j-1)X_6^\bullet (j)X_a^\bullet (j+1),  \quad &&a=8,6,   \\ 
&X_3^\bullet (j-2)X_8(j-1)X_4^\bullet (j-1)X_1^\bullet (j-1)X_6(j).
\quad &&                     
\endalignat$$
Here we used the multiplicative notation and the embedding
 $\rho_1\subset\pi$ is denoted with bullets.
\endproclaim

\proclaim{Lemma 13}
Let $\pi$ be one of the partitions of length 4 listed in Lemma 12. Set
$$
\pi_0=\cases X_a(i)X_3(j-1)X_5(j)X_1(j) \text{ \ if \ }
\pi=X_a(i)X_3(j-1)X_4(j)X_1(j),\\
X_a(i)X_8(j-1)X_5(j-1)X_6(j) \text{ \ if \ } \pi=X_a(i)X_8(j-1)X_4(j-1)X_6(j).
\endcases$$ 
Let $\pi'$ be such that $\wt (\pi')=\wt(\pi)$ and
$\pi_0\prec \pi'\prec \pi$. Then there is no embedding $\rho\subset\pi'$
such that $\rho$ is of the form 
$X_3(j'-1)X_5(j')X_1(j')$ or $X_8(j'-1)X_5(j'-1)X_6(j')$.
\endproclaim

\demo{Proof}
First note that $\sh (\pi')=\sh (\pi)$, that $\wt (X_a)$ is a root or zero,
and that 
$$
\wt \bigl (X_3(j'-1)X_4(j')X_1(j')\bigr )=2\alpha_2+\alpha_1, 
\wt \bigl (X_8(j'-1)X_(j'-1)X_6(j')\bigr )=-2\alpha_2-\alpha_1.
$$
Let $\pi=X_a(i)X_3(j-1)X_4(j)X_1(j)$. Then there is no root $\beta$ such that
$\wt (X_a)+2\alpha_2+\alpha_1=\beta-2\alpha_2-\alpha_1$, and hence $\pi'$
cannot
contain a partition of the form $X_8(j'-1)X_5(j'-1)X_6(j')$.

Assume that $\pi'=X_{a'}(i')X_3(j'-1)X_5(j')X_1(j')$. If $a\neq 4$ and
$a\neq 5$,
then $\wt (\pi')=\wt(\pi)$ implies that $a'=a$, 
and it is easy to see that this implies that $\pi'=\pi_0$.
If $a= 4$ or $a= 5$, then $\wt (\pi')=\wt(\pi)$ gives $a'=4$ or $a'=5$
and it is easy to see that then $\pi_0\prec \pi'\prec \pi$ does not hold.

The case when $\pi=X_a(i)X_8(j-1)X_4(j-1)X_6(j)$ is similar.\qed 
\enddemo

\proclaim{Lemma 14}
Let $\pi$ be one of the partitions listed in Lemma 12 and let
$\rho_1, \rho_2\subset\pi$. Then
$$
u(\rho_1 \subset \pi) \in \Bbb F^\times u(\rho_2 \subset \pi) +
\Cal R_{(\pi)}.\tag 8
$$
\endproclaim

\demo{Proof}
Let $\pi=X_3(j-1)X_5(j)X_4(j)X_1(j)$. {}From (7) we get
$$\aligned
&r\bigl(X_3(j-1)X_4(j)X_1(j)\bigr)X_5(j)\\
&=X_3(j-1)X_5(j)r\bigl(X_5(j)X_1(j)\bigr)
-r\bigl(X_3(j-1)X_5(j)\bigr)X_5(j)X_1(j)+\Cal R_{(\pi)}\\
&=r\bigl(X_3(j-1)X_5(j)\bigr)X_4(j)X_1(j)+\Cal R_{(\pi)},
\endaligned$$
and this implies (8) for any $\rho_1, \rho_2\subset\pi$.
 
Let $\pi=X_a(i)X_3(j-1)X_4(j)X_1(j)$ and assume that there is no embedding 
$\rho\subset\pi$ such that $\rho$ is of the form 
$X_3(j'-1)X_5(j')$ or $X_5(j')X_1(j')$, that is $\pi\neq
X_3(j-1)X_5(j)X_4(j)X_1(j)$.
Clearly there is no embedding 
$\rho\subset\pi$ such that $\rho$ is of the form 
$X_8(j'-1)X_5(j'-1)$ or $X_5(j'-1)X_6(j')$. Hence for every two
embeddings $\rho_1, \rho_2\subset\pi$, $\rho_1, \rho_2\in\lt (\bar R)$, the
relation (5) holds.
 
{}From the definition (2) we have
$$\aligned
 u&=u\bigl(X_3(j-1)X_4(j)X_1(j)\subset
X_a(i)X_3(j-1)X_4(j)X_1(j)\bigr)\\
&=X_a(i)r\bigl(X_3(j-1)X_4(j)X_1(j)\bigr)+v\\
&=X_a(i)X_3(j-1)r\bigl(X_5(j)X_1(j)\bigr)-X_a(i)r\bigl(X_3(j-1)X_5(j)\bigr)
X_1(j)+v,
\endaligned\tag 9
$$
where $v\in \Cal R_{(\pi)}$ may arise {}from commuting elements to the indicated
positions. Note that the leading term of $u$ is $\pi$. 

If $i=j-2$, then $a=3,1$, and by using the relation (5) for two embeddings
in $X_a(j-2)X_3(j-1)X_5(j)$ we get {}from (9)
$$
u\in X_a(i)X_3(j-1)r\bigl(X_5(j)X_1(j)\bigr)
-r\bigl(X_a(i)X_3(j-1)\bigr)X_5(j)X_1(j)+\Cal R_{(\pi_0)},
$$
where $\pi_0$ is defined in Lemma 13. By using two different expansions,
we get $u\in \Cal R_{(\pi_0)}$, i.e.
$$
u\in \sum_{\pi_0\prec \pi'\preccurlyeq \pi}\sum_{\rho\in\Cal E(\pi')}
c_{\rho,\pi'}\,u(\pi'/\rho) r(\rho)+\Cal R_{(\pi)}\tag 10
$$
for some coefficients $c_{\rho,\pi'}\in \Bbb F$. By Proposition 10 and Lemma 13
for all embeddings $\rho\subset\pi'$, $\pi_0\prec \pi'\preccurlyeq \pi$, we can
apply (5), so for each $\pi'$ we choose a particular $\rho_{\pi'}\subset\pi'$,
$\rho_{\pi'}\in \lt (\bar R)$, and
(by using (5) if necessary) we get
$$
u\in \sum_{\pi_0\prec \pi'\preccurlyeq \pi}
c_{\pi'}\,u(\pi'/\rho_{\pi'}) r(\rho_{\pi'})+\Cal R_{(\pi)}
$$
for some coefficients $c_{\pi'}\in \Bbb F$. Since $\lt (u)=\pi$, 
we get by induction (cf. [MP2, Lemma 9.4])
that $c_{\pi'}=0$ for $\pi_0\prec \pi'\prec \pi$, and hence
$$
X_a(j-2)r\bigl(X_3(j-1)X_4(j)X_1(j)\bigr)\in 
c_{\pi}\,u(\pi/\rho_{\pi}) r(\rho_{\pi})+\Cal R_{(\pi)}\tag 11
$$
for some particular embedding $\rho_{\pi}\subset\pi$, $\rho_{\pi}\in
\lt (\bar R)$.
Now (11) and (5) imply (8) for any $\rho_1, \rho_2\subset\pi$.

If $i=j-1$, then $a=7,5,4,3,2,1$, and by using the relation (5) for two
embeddings
in $X_a(j-1)X_3(j-1)X_5(j)$ we get {}from (9)
$$
u\in X_a(i)X_3(j-1)r\bigl(X_5(j)X_1(j)\bigr)
- r\bigl(X_a(i)X_3(j-1)\bigr)X_5(j)X_1(j)+\Cal R_{(\pi_0)}
$$
for all $a=7,5,4,3,2,1$.

If $i=j$, then $a\in \{8,\dots,1\}\backslash \{5\}$. 
By using the relation (5) for two embeddings
in $X_a(j-1)X_3(j-1)X_5(j)$ we get {}from (9)
$$
u\in X_a(i)X_3(j-1)r\bigl(X_5(j)X_1(j)\bigr)
-r\bigl(X_a(i)X_3(j-1)\bigr)X_5(j)X_1(j)+\Cal R_{(\pi_0)}
$$
for $a=8,7,6,3$.
By using (5) for two embeddings
in $X_5(j)X_a(j)X_1(j)$ we get {}from (9)
$$
u\in X_3(j-1)X_5(j)r\bigl(X_a(i)X_1(j)\bigr)
-r\bigl(X_3(j-1)X_5(j)\bigr)X_a(i)X_1(j)+\Cal R_{(\pi_0)}
$$
for $a=2,1$. For $a=4$, by using (5), we get {}from (9) that
$$
u= C X_3(j-1)r\bigl(X_5(j)X_a(i)\bigr)X_1(j)+\Cal R_{(\pi_0)}.
$$
Since $\pi_0\prec \pi=\lt(u)$, this implies that $C=0$ and 
$u\in\Cal R_{(\pi_0)}$.

If $i=j+1$, then $a\in \{8,\dots,1\}$, and
by using (5) for two embeddings in \newline
$X_5(j)X_1(j)X_a(j+1)$ we get {}from (9)
$$
u\in X_3(j-1)X_5(j)r\bigl(X_1(j)X_a(i)\bigr)
-r\bigl(X_3(j-1)X_5(j)\bigr)X_1(j)X_a(i)+\Cal R_{(\pi_0)}
$$
for all $a\in \{8,\dots,1\}$.

In either of these cases we get $u\in \Cal R_{(\pi_0)}$, i.e. (10), and we
argue 
as before that (8) holds for any $\rho_1, \rho_2\subset\pi$.

The case when $\pi=X_a(i)X_8(j-1)X_4(j-1)X_6(j)$ is similar.

Finally, let us consider the case 
$\pi=X_3(j-2)X_8(j-1)X_4(j-1)X_1(j-1)X_6(j)$,
$\rho_1=X_8(j-1)X_4(j-1)X_6(j)$. By applying relations for the embeddings
in $X_3(j-2)X_8(j-1)X_4(j-1)X_1(j-1)$ and $X_3(j-2)X_8(j-1)X_4(j-1)X_6(j)$
that we have already proved we get
$$\align
&r\bigl(X_3(j-2)X_4(j-1)X_1(j-1)\bigr)X_8(j-1)X_6(j)\\
&=r\bigl(X_3(j-2)X_8(j-1)\bigr)X_4(j-1)X_1(j-1)X_6(j)+\Cal R_{(\pi)}\\
&=X_3(j-2)X_1(j-1)r\bigl(X_8(j-1)X_4(j-1)X_6(j)\bigr)+\Cal R_{(\pi)},
\endalign$$
so (8) holds in this case as well.\qed
\enddemo

Lemma 14 is the last step in the proof of Theorem 11.

\demo{Proof of Theorem A}
Note first  that the order $\preccurlyeq$ on $\Cal P(\bar B_{<0})$ is
a (reverse) well order which behaves well with respect to
multiplications (cf. Lemmas 6.2.1 and 6.2.2 in \cite{MP2}).
Since the relations $r(\rho)$, $\rho\in\Cal R$, vanish on $L(\Lambda_0)$,
we see by induction that the set of vectors
$$
u(\pi)v_0,\qquad \pi\in\Cal P(\bar B_{<0})\backslash (\Cal P(\bar B)\Cal R)
$$
form a spanning set of $L(\Lambda_0)$. Of course, here
$\pi\in\Cal P(\bar B_{<0})\backslash (\Cal P(\bar B)\Cal R)$ means that
$u(\pi)$ satisfies the difference $\Cal R$ condition defined previously.

Since $L(\Lambda_0)=N(\Lambda_0)/N^1(\Lambda_0)$, linear independence
of this spanning set will follow {}from Proposition 6.3.2 in \cite{MP2},
provided we can construct a basis of $N^1(\Lambda_0)$ parametrized
by $\Cal P(\bar B_{<0})\cap (\Cal P(\bar B)\Cal R)$.

Since $N^1(\Lambda_0)=U(\tilde{\goth g}_{<0})\bar{R}\1 $ 
(see Section 5 in \cite{MP2}), the set of vectors
$$
u(\rho \subset \pi)\1,\qquad\rho \in \Cal R, \ 
\pi \in \Cal P(\bar B_{<0})\cap (\Cal P(\bar B)\Cal R), \ \rho \subset \pi
$$
is a spanning set of $N^1(\Lambda_0)$.
Now for each $\pi$ we choose 
precisely one $\rho=\rho(\pi)$ such that $\rho\subset \pi$.
By using Theorem 11, and induction on $\preccurlyeq$, we
see that
$$
u(\rho(\pi) \subset \pi)\1, \qquad \pi \in 
\Cal P(\bar B_{<0})\cap (\Cal P(\bar B)\Cal R)
$$
is a spanning set as well. Since 
$\lt_{N(\Lambda_0)}(u(\rho(\pi) \subset \pi)\1)=u(\pi)\1$, this spanning
set is obviously linearly independent, and hence the desired basis of
$N^1(\Lambda_0)$.
\qed
\enddemo


\head{One combinatorial identity }\endhead

Let $A$ be a nonempty set and denote by $\Cal P (A)$ the set
of all maps
$f : A \rightarrow \Bbb N$, where $f_a=f(a)$ equals zero for
all but finitely
many $a \in A$.
We say that $f$ is a partition and for $f_a > 0$ we
say that
$a$ is a part of $f$. 

For nonempty subsets $A_1, \dots ,A_s \subset \Bbb Z_{>0}$,
$s \ge 1$, let
$
 A = A_1 \bigsqcup \dots \bigsqcup A_s
$
be a disjoint union of sets. We call the elements of $\Cal P(A)$ colored
partitions
with parts in $A$, where for $c \in \{1, \dots ,s\}$
and $j \in A_c$ we say that
$j$ is of color $c$ and of degree $|j| = j \in
\Bbb Z_{>0}$. We define the degree $|f|$ of $f \in \Cal P(A)$ as
$|f| = \sum_{a
\in A}\,|a|f_a$ and we say that
$f$ is a colored partition in $s$ colors  of the nonnegative
integer $|f|$.

The explicit construction of the basic module for the affine Lie
algebra $\goth{sl}(3,\Bbb C)\sptilde$, given by Theorem A, implies the
following
combinatorial identity of
Rogers-Rama\-nu\-jan type:

\proclaim{Theorem B}  The number of partitions $f$ in one
color (say $A=\Bbb Z_{>0}$) of a
nonnegative integer $n$ such that each part appears at most
twice (i.e $f_a\le 2$) equals 
the number of partitions $f$ in three colors 
(say $A=\{r,\u r,\u{\u r} \mid r\in \Bbb Z_{>0}\}$) of $n$ 
such that each part appears at most once (i.e. $f_a\le 1$), 
but subject to the conditions

$$\gather
f_{\u{\u r}} >0 \quad\text{implies}\quad |\,\u{\u r}\,|
\equiv \pm 1 \mod{3}, 
\tag 1\\
|r|>1,\quad |\u{r}|>1, \quad |\,\u{\u r}\,| \ne 2 ,\tag 2\\
f_a+f_b\le 1 \quad \text{for}\quad  \big ||a|- |b|\big | \le 1,  
\tag 3
\endgather$$
and the conditions
$$\gathered
f_{3i+2}+f_{3i+3}+f_{\u{3i+4}}+f_{3i+5}\le 1,\\
f_{\u{3i+2}}+f_{\u{3i+3}}+f_{\u{\u{3i+4}}}+f_{\u{3i+5}}\le
1,\\
f_{\u{\u{3i+2}}}+f_{\u{3i+3}}+f_{3i+4}+f_{\u{\u{3i+5}}}\le
1,\\
f_{3i+1}+f_{\u{3i+2}}+f_{3i+3}+f_{3i+4}\le 1,\\
f_{\u{3i+1}}+f_{\u{\u{3i+2}}}+f_{\u{3i+3}}+f_{\u{3i+4}}\le
1,\\
f_{\u{\u{3i+1}}}+f_{3i+2}+f_{\u{3i+3}}+f_{\u{\u{3i+4}}}\le
1,
\endgathered \tag 4
$$
$$
\gathered
f_{\u{3i+1}}+f_{3i+3}+f_{\u{\u{3i+5}}}\le 2,\\
f_{\u{\u{3i+1}}}+f_{3i+3}+f_{\u{3i+5}}\le 2
\endgathered \tag 5
$$
for all $i$ in $\Bbb N$.
\endproclaim
\demo{Proof}
By using the Lepowsky numerator formula for the
 Weyl-Kac character formula (cf. \cite{K}, \cite{L}), we can write the
principally
specialized character of the basic $\goth{sl}(3,\Bbb C)\sptilde$-module
as the infinite product
$$
 \prod_{r\not\equiv 0\mod{3}}\, (1 - q^r)^{-1}=
\prod_{r\geq 1}\, (1 + q^r+q^{2r}).
$$
This product can be interpreted as the generating function
$\sum_{n\geq 0} c_nq^n$
of the partition function
$c_n$ counting the number of partitions $f$ of a nonnegative integer $n$ 
such that each part $r$ appears at most twice.

On the other hand, using Theorem A we can describe the principally
specialized character
of the basic module $L(\Lambda_0)$: Consider the isomorphism of
monoids $\varphi \: \Cal P(\bar B_{<0})\to \Cal P(A)$ defined by the bijection
$$
\varphi \: \bar B_{<0} \ \longrightarrow 
A = \{j\mid j\geq 2\} \bigsqcup
\{\u j \mid j\geq 2\} \bigsqcup 
\{\u{\u{3i\pm 2}} \mid i\geq 1\},
$$
$$\alignat2
&X_1(-i)\mapsto \u{\u{3i-2}}, \qquad\qquad &&X_8(-i)\mapsto \u{\u{3i+2}},\\
&X_2(-i)\mapsto {{3i-1}}, \qquad\qquad &&X_7(-i)\mapsto {{3i+1}},\\
&X_3(-i)\mapsto {\u{3i-1}}, \qquad\qquad &&X_6(-i)\mapsto {\u{3i+1}},\\
&X_4(-i)\mapsto {{3i}}, \qquad\qquad &&X_5(-i)\mapsto {\u{3i}}.
\endalignat$$
In particular, we (would) have $f_1=X_7(0)\mapsto 1$, $f_2=X_6(0)\mapsto
\u{1}$, 
$f_0=X_1(-1)\mapsto \u{\u{1}}$, and by setting 
$|{{n}}|=|{\u{n}}|=|\u{\u{n}}|=n\in\Bbb Z_{>0}$ the map
$$
X_s(-i)\mapsto \deg X_s(-i)=|\varphi (X_s(-i))|
$$
induces the principal specialization of the character of the basic
module: the degree of the vector $X_{i_1}(j_1)X_{i_2}(j_2)\dots
X_{i_s}(j_s)v_0$
is the sum $\sum \deg X_{i_r}(j_r)$.

Note that $\Cal P(\bar B_{<0})\backslash (\Cal P(\bar B)\Cal R)$ is a
partition ideal (cf. \cite{A}, \cite{MP2, 11.2}) defined by the difference
$\Cal R$ conditions.
The condition
that a colored partition $\pi=(\pi_a\mid a\in\bar B_{<0})$ does not
contain $X_s(j)X_s(j)$ for $j<0$ and $s=1,\dots, 8$ can be written as
$\pi_{X_s(j)}\leq 1$, meaning that the part $X_s(j)$ may appear in $\pi$
at most once. Since $\varphi$ is an isomorphism, for the corresponding
colored partition $f=\varphi(\pi)= \pi\circ\varphi^{-1}$, $f=(f_a\mid a\in A)$,
this condition reads
that the part $a$ may appear in $f$ at most once. Now the condition that
$\pi$ does not contain $X_1(-i-1)X_1(-i)$ can be written as
$\pi_{X_1(-i-1)}+\pi_{X_1(-i)}\leq 1$, and the corresponding condition for
$f$ as $f_{\u{\u{3i+1}}}+f_{\u{\u{3i-2}}} \leq 1$. 

If we write down all difference $\Cal R$ conditions for $\pi$ and
the corresponding conditions for $f$, we can arrange them to take the form
(3)--(5), as stated in the theorem. The conditions (5) come {}from the two
cubic terms in $\Cal R$, the conditions (2) represent the ``initial''
conditions $f_1\cdot v_0=0$ and $f_2\cdot v_0=0$ for the basic module.
\qed
\enddemo
\remark{Remark}
Let us remark that the conditions (3) and (4) stated in Theorem B represent
over fifty conditions of the form $f_a+f_b\leq 1$. We tried to write them
in a compact and comprehensive form, but still without
understanding in what form this identity might be generalized to
higher levels or higher ranks.
\endremark


\vfill\eject


\Refs

\widestnumber\key{AAG}

\ref\key{AAG}
\by K.~Alladi, G.~E.~Andrews and B.~Gordon
\paper Refinements and generalizations of Capparelli's
conjecture on partitions
\jour J. Algebra
\vol 174
\yr 1995
\pages 636-658
\endref

\ref\key{A}
\by G. E. Andrews
\book The theory of partitions, {\rm Encyclopedia of math.
and appl.}
\publ Addison-Wesley
\publaddr Amsterdam
\yr 1976
\endref


\ref\key{C}
\by S. Capparelli
\paper A construction of the level 3 modules for the affine Lie
algebra $A_2^{(2)}$ and a new combinatorial identity of the
Rogers-Ramanujan type
\jour Trans. AMS
\vol 348
\yr 1996
\pages 481--501
\endref

\ref\key{FLM}
\by I. B. Frenkel, J. Lepowsky and A. Meurman
\book Vertex operator algebras and the Monster,
{\rm Pure and Applied Math.}
\publ Academic Press
\publaddr San Diego
\yr 1988
\endref

\ref\key{K}
\by V. G. Kac
\book Infinite-dimensional Lie algebras {\rm 3rd ed.}
\publ Cambridge Univ. Press
\publaddr Cambridge
\yr 1990
\endref

\ref\key{L}
\by J. Lepowsky
\paper Application of the numerator formula to $k$-rowed plane partitions
\jour Advances in Math.
\vol 35
\yr 1980
\pages 179--194
\endref

\ref\key{LP}
\by J. Lepowsky and M. Primc
\paper Structure of the standard modules for the affine Lie algebra $A_1^{(1)}$
\jour Contemporary Math.
\vol 46
\yr 1985
\endref

\ref\key{LW}
\by J. Lepowsky, R. L. Wilson
\paper The structure of standard modules, I: Universal algebras and
the Rogers-Ramanujan identities
\jour Invent. Math.
\vol 77
\yr 1984
\pages 199--290
\moreref
\paper  \ II: The case $A_1^{(1)}$, principal gradation
\jour Invent. Math.
\vol 79
\yr 1985
\pages 417--442
\endref

\ref\key{Ma}
\by M. Mandia
\paper Structure of the level one standard modules for the affine Lie algebras
$B_l^{(1)}$, $F_4^{(1)}$ and $G_2^{(1)}$
\jour Memoirs American Math. Soc. 
\vol 362
\yr 1987
\endref

\ref\key{MP1}
\by A. Meurman, M. Primc
\paper Annihilating ideals of standard modules of $\frak{sl}(2,\C)\sptilde$
and combinatorial identities
\jour Advances in Math.
\vol 64
\yr 1987
\pages 177--240
\endref

\ref\key{MP2}
\by A. Meurman and M. Primc
\paper Annihilating fields of standard modules of $\frak{sl}(2,\Bbb C)\sptilde$
and combinatorial identities
\jour Memoirs American Math. Soc.
\vol 652
\yr 1999
\endref

\ref\key{Mi}
\by K. C. Misra
\paper Level one standard modules for affine symplectic Lie algebras
\jour Math. Ann. 
\vol 287
\yr 1990
\pages 287--302
\endref

\endRefs
\enddocument